\documentclass[a4paper,12pt]{amsart}
\usepackage{amsmath, amssymb, amsfonts,
 epsfig, mathrsfs, enumerate, ifthen,
 graphicx,
 tikz}

 \usepackage{tabularx}
\usepackage{multirow}
\nonstopmode \numberwithin{equation}{section}
\usepackage{lineno,hyperref,graphics,mathrsfs}
\usepackage{tkz-euclide}
\newtheorem{theorem}{Theorem}
\newtheorem{lemma}{Lemma}
\newtheorem{problem}{Problem}
\newtheorem{definition}{Definition}
\begin{document}
\title{On determining number and metric dimension of zero-divisor graphs}
\author{Muhammed Sabeel K}
\address{Muhammed Sabeel K, Department of Mathematics, National Institute of Technology Calicut, Kozhikode~\textnormal{673601}, India. }
\email{sabeel.math@gmail.com}
	
\author{Krishnan Paramasivam}
\address{Krishnan Paramasivam, Department of Mathematics, National Institute of Technology Calicut,  Kozhikode~\textnormal{673601}, India. }
\email{sivam@nitc.ac.in}
%\begin{abstract}
%	In this paper we express  zero-divisor graph of $\mathbb{Z}_n$ as generelized join of a set of graph and use this to  calculate the deteriming number and metric dimension of zero-divisor graph of $\mathbb{Z}_n$. 
%\end{abstract}
\maketitle
\pagestyle{myheadings}
\markboth{Muhammed Sabeel K and Krishnan Paramasivam }{ }
\begin{abstract}
	In this article, explicit formulas for finding the determining number and the metric dimension of the zero-divisor graph of $\mathbb{Z}_n$ and non-Boolean semisimple rings are given.  In  the case of Boolean rings, an upper bound of the determining number and the metric dimension of  zero-divisor graph is  determined.  Further, the determining number and the metric dimension of some important graphs other than zero-divisor graphs, are proved and the open problem by Boutin \cite{boutin}, regarding the determining number of graphs is settled.
\end{abstract}
\vskip .3cm \noindent In 2006, the concept of determining number of a graph was introduced and defined by Boutin \cite{boutin}, and in the same year, an equivalent graph parameter, namely  the fixing number of a graph was studied independently by Erwin and Harary \cite{harr}.  In 1975, the concept of metric dimension was defined by Slater \cite{slater} and independently studied by Harary and Melter\cite{hararymelter} in 1976. It was  Harary and Melter \cite{hararyandmel} coined the name metric dimension for the same. The metric dimension is also called the distance-determining number as well as the locating number. In 1988, The zero-divisor graph of a commutative ring was introduced by Beck \cite{beck1988}, and later in 1999, a modified definition was given by Anderson and Livingston \cite{anderson1999zero}. Throughout this article, the definition of zero-divisor graph by Anderson and Livingston \cite{anderson1999zero} is adopted. 
\vskip .3cm \noindent In this article, we emphasize studying two graph parameters namely, determining number and metric dimension of zero-divisor graph of some classes of commutative rings. The study of metric dimension in zero-divisor graphs was already carried out by Pirzada et al. \cite{pirzadametric} and lots of results are proved in this direction. Apart from the usual approach,  we are  looking at the problem of finding metric dimension via  its relationship with determining number. For this, we represent the zero-divisor graph of different classes of graphs as a generalized join of certain graphs and utilize the structure to identify the determining number and hence to find the metric dimension. 
 In \cite{pirzadametric}, Pirzada et al. determined the metric dimension of $\Gamma(\mathbb{Z}_n)$ when   $n=2p$ or $n=p^k$ for a prime $p$. In this article, the determining number of $\Gamma(\mathbb{Z}_n)$ is found for all integers $n$. In section 2, we discuss the same for semisimple commutative rings and  obtain certain interesting results regarding the finiteness of determining number of zero-divisor graphs of an infinite ring. In section 4, we prove some miscellaneous results regarding the determining number of graphs irrespective of zero-divisor graph.\\
\section{Definition and Terminology}
In this article, we consider simple and undirected graphs. A graph $\Gamma$ is an ordered pair $(V(\Gamma), E(\Gamma))$, where the vertex set $V(\Gamma)$ is a non-empty set and the edge set $E(\Gamma)$ is a subset of the set of all two-elements subsets of $V(\Gamma)$. A graph is finite if its vertex set is finite and a graph is infinite graph if its vertex set is not finite. For any vertex $v$ of a graph  $\Gamma$, $N_{\Gamma}(v)$ or shortly, $N(v)$ is the neighborhood of $v$ in $\Gamma$, which is the set of all vertices adjacent to $v$ in $\Gamma$, and $deg_{\Gamma}(v)$ or shortly, $deg(v)$ is the degree of $v$ in $\Gamma$ and is equal to $|N_{\Gamma}(v)|$. For any two vertices $u$ and $v$ in $\Gamma$, the distance between $u$ and $v$ is the length of the shortest path connecting $u$ and $v$. A graph is complete if any two vertices are adjacent and a graph is an empty graph if no two vertices are adjacent. A complete graph on $n$ vertices is denoted by $K_n$ whereas an empty graph on $n$ vertices is denoted by $\overline{K}_n$. If $S$ is a subset of the vertex set of $\Gamma$, then the subgraph induced by $S$, denoted by $\langle S \rangle $, is the subgraph of $\Gamma$ with vertex set $S$ and two vertices $u$ and $v$ are adjacent in  $\langle S\rangle $ if and only if $u$ and $v$ are adjacent in $\Gamma$. A subset $S$ of vertices  of $\Gamma$ is a clique of $\Gamma$ if the subgraph $\langle S\rangle$ induced by $S$ is a complete graph and $S$ is an independent set if the subgraph $\langle S\rangle$ of $\Gamma$ induced by $S$ is an empty graph.\\
An automorphism of a graph $\Gamma$ is a bijection $\sigma$ from the vertex set of  $\Gamma$ to itself in such a way that
$u$ is adjacent to $v$ in $\Gamma$ if and only if $\sigma(u)$ is adjacent to $\sigma(v)$ in $\Gamma$. For any graph $\Gamma$, $Aut(\Gamma)$ stands for the group of all automorphisms of $\Gamma$.
For a graph $\Gamma$, a non-empty subset $\mathcal{O}$ of $V(\Gamma)$ is said to be a vertex orbit of $\Gamma$ if  $\mathcal{O}=\{\sigma(u): \sigma \in Aut(\Gamma)\}$ for some $u$ in $V(\Gamma)$.  A graph $\Gamma$ is vertex-transitive if $V(\Gamma)$ is the only vertex orbit of $\Gamma$. The stabilizer of a vertex $v$ of a graph $\Gamma$, is $stab(v)=\{\sigma \in Aut(\Gamma):\sigma(v) =v\}$. The vertex stabilizer of a set of vertices $S$ of a graph $\Gamma$ is  $stab(S) =\{\sigma \in Aut(\Gamma):\sigma(v) =v, \forall v \in S\}$.\\

\par For a commutative ring $R$ with 1, $\mathcal{Z}(R)$ denotes the set of all non-zero zero-divisors of $R$ and $Ann(R)$ denotes the set of non-zero ideals with non-zero annihilator. Note that in this article, the zero-element $0$ is not included in $\mathcal{Z}(R)$, and the zero-ideal $\langle\bf{0}\rangle$ is not included in $Ann(R)$. A ring $R$ is Boolean if $x^2=x$ for any $x\in R$. A commutative ring is Noetherian if it satisfies ascending chain condition of ideals. A prime ideal $\mathfrak{p}$ of $R$ is an associated prime if  $\mathfrak{p}=ann(x)$  for some  $x\in R$. The set of all associated primes ideals in $R$ is denoted $Ass(R)$. In this paper, we only consider commutative rings, therefore by a ring we always mean a commutative ring until otherwise it is mentioned as non-commutative.\\
For an element $a$ of a group $G$, $\langle a \rangle$ denotes the subgroup generated by $a$. For any non-empty set $X$, $S_X$ denotes the group of all permutations of $X$, and $S_n$ denotes the symmetric group on $n$ symbols. For any finite set $X$, $|X|$ denotes the cardinality of $X$ and $X^c$ represents the complement of the set $X$ with respect to the appropriate superset of $X$. We use $[1,k]$ to denote the set $\{1,2,\cdots,k\}$.  For any real number $\alpha$, $\lfloor \alpha\rfloor$ denotes the greatest integer  function of $\alpha$.

\begin{definition}\cite{boutin} Let $\Gamma=(V(\Gamma),E(\Gamma))$ be a graph.	A subset $S \subseteq V(\Gamma)$ is said to be a determining set of $\Gamma$ if whenever $\sigma_1, \sigma_2 \in Aut(\Gamma)$  such that $\sigma_1(u) = \sigma_2(u)$ for all $u\in S$, then $\sigma_1(v) = \sigma_2(v)$ for all $v\in V (\Gamma)$. The determining number of a graph $\Gamma$ is defined as  $$Det(\Gamma)=\min_{S \subset V(\Gamma)}\{|S| : S \textnormal{~is~ a ~determining~ set~ of~} \Gamma\}.$$
\end{definition}
 
For any graph $\Gamma$, it is obvious that $V(\Gamma)$ is a determining set of $\Gamma$. Also any set of vertices $S$ with $|S|=|V(\Gamma)|-1$ is a determining set of $\Gamma$. Therefore $0\le Det(\Gamma)\le |V(\Gamma)|-1$.  

The notion of fixing sets and fixing number are analogous to the concept of determining sets and determining number. Below  we shall see that both these concepts of fixing sets and determining sets are equivalent.

\begin{definition}\cite{harr}\label{fix}
    A vertex $v$ of a graph $\Gamma$ is said to be fixed by an automorphism $\sigma \in Aut(\Gamma)$ if $\sigma \in stab(v)$. A set of vertices of $\Gamma$ is a fixing set of $\Gamma$ if $stab(S)$ is trivial, and we say that $S$ fixes $\Gamma$. The fixing number  of $\Gamma$ is $$\textnormal{fix}(\Gamma)=\min_{S \subset V(\Gamma)}\{|S| : S \textnormal{~is~ a ~fixing~ set~ of~} \Gamma\}.$$
\end{definition}

Gibbons and Laison \cite{fixgibbon} proved that the determining sets and the fixing sets are not different; but they are equivalent.

\begin{lemma}\cite{fixgibbon}\label{fixdet}
    A set of vertices is a fixing set if and only if it is a determining set
\end{lemma}
We later focus on the concept of metric dimension, which is equivalent notion of the resolving set of a graph. 
\begin{definition}\cite{slater}
 Let $\Gamma=(V(\Gamma),E(\Gamma))$ be a graph and $W =\{w_1,\cdots, w_k\}$ be an ordered subset  of $V(\Gamma)$. Then for a vertex $v$ of $\Gamma$, the metric representation of $v$ with respect to $W$ is the $k$-vector $D_W(v) = (\zeta_1, \zeta_2, 
 \cdots,\zeta_k)$, where $\zeta_i$ is the distance between  $v$ and $w_i$ in $\Gamma$.
\end{definition}
\begin{definition}\cite{slater}
 Let $\Gamma$ be a graph. A set of vertices $S$ is said to  resolving set or distance-determining set of  $\Gamma$ if every vertex of $\Gamma$ is uniquely
determined by the metric representation of the vertices with respect to $S$, that is, $D_S(u) \ne  D_S(v) $ for all $u, v\in  V(\Gamma)$ with $u\ne  v$, where $S$ assumed some order. The metric dimension of $\Gamma$  is the minimum cardinality of a resolving set of $\Gamma$,  denoted by $dim_M(\Gamma)$ 
\end{definition}
\begin{lemma}\label{detmd}\cite{detandmd}
If $S \subseteq V(\Gamma)$ is a resolving set of  a graph $\Gamma$, then $S$ is a determining set of $\Gamma$. In particular, $Det(\Gamma) \le dim_M(\Gamma)$.
\end{lemma}
We first discuss the determining number and the metric dimension of certain well-known classes of graphs. For the complete graph $K_n$, any set $S$ with $|S|=n-1$ is a determining set of $K_n$, we have $Det(K_n)=n-1=dim_M(K_n)$. For a path $P_n$, a single pendant vertex of $P_n$ forms a determining set of $P_n$ and therefore $Det(P_n)=1=dim_M(P_n)$. For a cycle $C_n$, $Det(C_n)=2=dim_M(C_n)$, because any pair of two vertices in $C_n$ form a determining set and a resolving set. \\

The following graph operation namely, the generalized join of graphs by Sabidussi \cite{sabidussi1961graph}, is often used in this article.
\begin{definition}\cite{sabidussi1961graph} \label{genjoin} Let $\Gamma$ be a given graph and $\{\Lambda_{\alpha}\}_{\alpha \in V(\Gamma)} $ be a collection of graphs indexed $V(\Gamma)$. Then the generalized join of $\Gamma$ with $\{\Lambda_{\alpha}\}_{\alpha \in V(\Gamma)} $ is a graph  $\tilde{\Gamma}$  with vertex set $V(\tilde{\Gamma})=\{(x,y):x\in V(\Gamma) ~\textnormal{and}~ y\in V(\Lambda_x) \}$ and two vertices $(x,y)$ and $(x',y')$ are adjacent if and only if either $x$ is adjacent to $x'$ in $E(\Gamma)$ and $y$ is adjacent to $y'$ in $E(\Lambda_x)$ or $x=x'$ and $y$ is adjacent to $y'$ in $E(\Lambda_x)$.  If $\Gamma$ has $m$ vertices, then $\Gamma$-join of the collection $\{\Lambda_1,\Lambda_2, \cdot \cdot\cdot,\Lambda_m\}$ is denoted by  $\Gamma[\Lambda_1,\Lambda_2,\cdots,\Lambda_m]$.
\end{definition}
The following definition of zero-divisor graph by Anderson and Livingston \cite{anderson1999zero} is adopted in this article.

\begin{definition}\cite{anderson1999zero} Let $ R $ be a  ring. The zero-divisor graph $\Gamma(R)$, of $R$, is a simple graph having the vertex set $\mathcal{Z}(R)$ and two distinct elements, $x$ and $y$ of $\mathcal{Z}(R)$ are adjacent if and only if $xy=0$. \end{definition}

The following result determines the finiteness of the order of the zero-divisor graph.

\begin{lemma}\label{Ganes}\cite{ganesan}
For any  ring $R$. $\mathcal{Z}(R)$ is finite if and only if either $R$ is finite or is an integral domain. 
\end{lemma}

Define an equivalence relation $\lq \sim$'   on $\mathcal{Z}(R)$ in such a way that for any $x,y \in \mathcal{Z}(R)$, we define $x\sim y$  if and only if $ann(x)=ann(y)$. With respect to this equivalence relation $\lq \sim$', if $R_E$ denotes the set of all equivalence classes $[x]$ in $\mathcal{Z}(R)$, then $R_E$ forms a partition of $\mathcal{Z}(R)$.\\
\begin{definition}\cite{spiroff2011zero}
	For a  ring, $R$ the \textit{compressed zero-divisor graph} of $\Gamma_E(R)$ of $R$ is a simple graph with the vertex set $R_E$ and two distinct vertices $[x]$ and $[y]$ of $R_E$ are adjacent if and only if $xy=0$. This equivalence relation is the base of compressed zero-divisor graph.
\end{definition}

Another interesting graph structure associated with  a ring is the annihilating ideal graph, which is defined as follows.
\begin{definition}\cite{behboodi2011annihilating}
	Let $R$ be a  ring and $Ann(R)$ be the set of ideals with a non-zero annihilator. The \textit{annihilating-ideal graph} $\Gamma_{Ann}(R)$ of $R$ is the  graph with the vertex set $Ann(R)$, and two distinct vertices $I$ and $J$ are adjacent if and only if $IJ=\langle\bf{0}\rangle$.
\end{definition}

There are several structural properties of these three graphs zero-divisor graph, compressed zero-divisor graph, and annihilating ideal graphs, identified by different authors. To review further properties of these graphs, we refer \cite{anderson1999zero,spiroff2011zero,behboodi2011annihilating,newsurvayzdg}.

\section{Determining number and metric dimension of $\Gamma(\mathbb{Z}_n)$}
In this article, we analyze certain structural properties of zero-divisor graph of $\mathbb{Z}_n$ that helps to predict the determining number and the metric dimension of $\Gamma(\mathbb{Z}_n)$.\\
While discussing  the  zero-divisor graph of $\mathbb{Z}_n$, $n$ represents a positive integer other than 1 or a prime. In addition, the following notations and terminologies are used in this article.\\
For a positive integer $n>1$, $\phi(n)$  denotes the Euler totient function of $n$. We call $d$ is a proper divisor of $n$ if $d\ne 1$, $d\ne n$ and $d$ divides $n$. We use $\Upsilon(n)$ to denote the set of all divisors of $n$ and $\tau(n)$ denotes the number of proper divisors of $n$. Note that $1$ and $n$ are not included as the elements of $\Upsilon(n)$.  We further assume that the set of all integers modulo $n$ is
denoted by $\mathbb{Z}_n=\{{\bf 0}, {\bf 1},\cdots,{\bf n-1}\}$. Note that any element $\textbf{k}$ of $\mathbb{Z}_n$ is denoted by a bold letter and the corresponding integer $k$ with $0\le k<n$, is an element of $\mathbb{Z}$.
\begin{definition}
Let $n>1$ be any positive integer and let $d$ be any divisor of $n$ such that $1<d<n$. We define  
${\bf \Omega}_d$ to be the set of all ${\bf x}$ in $\mathbb{Z}_n$ such that greatest common divisor of the corresponding integer $x$ and $n$ is $d$. That is, $${\bf \Omega}_d=\{ {\bf x}\in \mathbb{Z}_n: gcd(x,n)=d\}.$$
\end{definition}
\noindent {\bf Note 2.1} Since $d \in {\Omega}_d$ for any divisor $d$ of $n$, the set ${\bf \Omega}_d$ is always non-empty.  Also, for any divisor $d$ of $n$, ${\bf \Omega}_d$ is a subset of $ \mathcal{Z}(\mathbb{Z}_n)$ and moreover ${\bf \Omega}_{d_1} \cup {\bf \Omega}_{d_2}\cup \cdots \cup {\bf \Omega}_{d_{\tau{(n)}}}= \mathcal{Z}(\mathbb{Z}_n)$. \\ \noindent {\bf Note 2.2} If  $d$ and $d'$ are any two distinct divisors of $n$, then ${\bf \Omega}_d\cap{\bf \Omega}_{d'}=\{\}$. In addition, if    $\mathcal{V}=\{{\bf \Omega}_d$: $d\in \Upsilon(n)\}$, then  one can see that $\mathcal{V}$ forms a partition of $\mathcal{Z}(\mathbb{Z}_n)$.\\
\noindent {\bf Note 2.3}. For any divisor $d$ of $n$, it can be verified that ${\bf \Omega}_d=\{ {\bf d}{\bf u}: {\bf u} \textnormal{~is a unit in~} \mathbb{Z}_n\}$.

%\begin{proposition}
%$|\Omega_d|=\phi(\frac{n}{d})$ 
%\end{proposition}
%\begin{proof}
%We have additive order of an element $x\in\mathbb{Z}_n$ will be, $O(x)=\frac{n}{gcd(x,n)}=\frac{n}{d}$. Therefore, $\Omega_d$ is collection of all elements of order $\frac{n}{d}$ in $\mathbb{Z}_n$. Therefore, $|\Omega_d|=\phi(\frac{n}{d})$.
%\end{proof}

%We can see that, for a fixed $d\in \Omega_d$, if  $x,y\in \Omega_d$ , then  for any $a\in R$, $ax=0$ if and only if $ay=0$.  Therefore, $ann(x)=ann(y)$ for $x,y \in \Omega_d$. Consequently $x\sim y$ if and only $x,y \in \Omega_d$ as a result $\Omega_d's$ are the equivalence classes corresponding to the relation $\sim$ 

From the following lemmas, one can realize  the sets ${\bf \Omega}_d$'s  guarantee to explore  certain important structural properties of the zero-divisor graph of $\mathbb{Z}_n$. First, we prove that for any $d\in \Upsilon(n)$ the set ${\bf \Omega}_d$ is a clique or an independent set in $\Gamma(\mathbb{Z}_n)$.
 
\noindent \begin{lemma}\label{3lemmas}
Let $d$ be a divisor of $n$. 
\begin{enumerate}
 \item[$(i)$] \label{CIset}  If $n|d^2$ then ${\bf \Omega}_d$ is  a clique $\Gamma(\mathbb{Z}_n)$ and ${\bf \Omega}_d$ is an independent set of $\Gamma(\mathbb{Z}_n)$ otherwise.
 \item[$(ii)$] \label{Vdsize}  The number of elements in ${\bf \Omega}_d$ is $\phi(\frac{n}{d})$.
 \item[$(iii)$] \label{degree}
	If $\textbf{x}\in {\bf \Omega}_d$, then  
	
\[ deg_{\Gamma(\mathbb{Z}_n)}({\bf x})= \begin{cases} 

      d-2 &  n\mid d^2.\\
      d-1 &  n\nmid d^2 .\\
   \end{cases}
\]
 \end{enumerate}
\end{lemma} 
%%%%%%%%%%%%%%%%%%%%%%%%%%%% 
\begin{proof}
$(i)$ Suppose that $\textbf{x},\textbf{y}\in {\bf \Omega} _d$ then there exist two units $\textbf{u}$ and $\textbf{v}$ of $\mathbb{Z}_n$ such that $\textbf{x}=\textbf{d}\textbf{u}$ and $\textbf{y}=\textbf{d}\textbf{v}$. Now $\textbf{x}\textbf{y}=\textbf{d}^2\textbf{u}\textbf{v}$, where $\textbf{u}\textbf{v}$ is a unit in $\mathbb{Z}_n$. Thus $\textbf{x}$ and $\textbf{y}$ are adjacent in $\Gamma(\mathbb{Z}_n)$ if and only if $n|d^2$.  Since $\textbf{x}$ and $\textbf{y}$ are arbitrary, the result follows.\\
 \noindent $(ii)$ One can see that $\Omega_d$ is a collection of all generators of the additive subgroup $\langle \textbf{d} \rangle$ of $\mathbb{Z}_n$ and hence $|{\bf \Omega}_d|$ is the number of generators of  $\langle \textbf{d} \rangle$. But $\langle \textbf{d} \rangle$ is isomorphic to  $\mathbb{Z}_{\frac{n}{d}}$. Thus $|{\bf \Omega}_d|=\phi(\frac{n}{d})$.

\noindent $(iii)$ If $\textbf{x}\in {\bf \Omega}_d$, then $\textbf{x}=\textbf{d}\textbf{u}$ for some unit $\textbf{u}$ in $\mathbb{Z}_n$. Now if $\textbf{y}\in \Omega_m$, where $m\in \Upsilon(n)$ and $m$ is a multiple of $\frac{n}{d}$, then $\textbf{x}\textbf{y}=\textbf{0}$. Suppose that $\textbf{x}\textbf{w}=\textbf{0}$ for some $\textbf{w} \in {\bf \Omega}_k$, where $k\in \Upsilon(n)$. Then $\textbf{d}\textbf{k}=\textbf{0}$ and hence $k$ is a multiple of $\frac{n}{d}$. \\
	If $n\nmid d^2$, then using the Lemma \ref{Vdsize}(ii) $$deg_{\Gamma(\mathbb{Z}_n)}(\textbf{x})=\sum_{ \substack{\frac{n}{d}|m\\ m\in \Upsilon(n)}}|{\bf \Omega}_m|=\sum_{\substack{\frac{n}{d}|m\\ m\in \Upsilon(n)}}\phi\bigl(\frac{n}{m}\bigr)=\sum_{\substack{i|d\\ i\ne 1}}\phi(i)=d-1.$$
	
	If $n| d^2$,  then $d$ is a multiple of $\frac{n}{d}$. Therefore when the counting is done as above $\textbf{x}$ will be included in  the $d-1$ number of adjacent vertices of the vertex $\textbf{x}$.  Since the zero-divisor  graph does not  admit  any loop, $deg_{\Gamma(\mathbb{Z}_n)}(\textbf{x})=d-2$, 
\end{proof}
\begin{lemma}
Let ${\bf x}, {\bf y}$ be any two vertices in $\Gamma(\mathbb{Z}_n)$. Then $deg_{\Gamma(\mathbb{Z}_n)}({\bf x})=deg_{\Gamma(\mathbb{Z}_n)}({\bf y})$ if and only if ${\bf x}, {\bf y} \in {\bf \Omega}_d$ for some $d\in \Upsilon(n)$.  
\end{lemma}
\begin{proof}
	Let $\textbf{x},\textbf{y} \in {\bf \Omega}_d$. Then $\textbf{x}=\textbf{d}\textbf{u}$ and $\textbf{y}=\textbf{d}\textbf{v}$, for some units $\textbf{u}$ and $\textbf{v}$ in $\mathbb{Z}_n$. Therefore, for any $\textbf{w}\in \mathcal{Z}(\mathbb{Z}_n)$, $\textbf{x}\textbf{w}=\textbf{0}$ if and only if  $\textbf{d}\textbf{w}=\textbf{0}$ if and only if  $\textbf{y}\textbf{w}=\textbf{0}$. Therefore,  $ann(\textbf{x})=ann(\textbf{y})$. By Lemma \ref{3lemmas}$(i)$, if $n|d^2$, then $deg_{\Gamma(\mathbb{Z}_n)}(\textbf{x})=|ann(\textbf{x})\setminus\{\textbf{0},\textbf{x},\textbf{y}\}|=|ann(\textbf{y})\setminus\{\textbf{0},\textbf{x},\textbf{y}\}|=deg_{\Gamma(\mathbb{Z}_n)}(\textbf{y})$ and if $n\nmid d^2$, then $deg_{\Gamma(\mathbb{Z}_n)}(\textbf{x})=|ann(\textbf{x})|-2=|ann(\textbf{y})|-2
	=deg_{\Gamma(\mathbb{Z}_n)}(\textbf{y})$. \\
	Conversely, if $deg_{\Gamma(\mathbb{Z}_n)}(\textbf{x})=deg_{\Gamma(\mathbb{Z}_n)}(\textbf{y})$ for $\textbf{x}\in \Omega_{d_1} $ and $\textbf{y}\in \Omega_{d_2}$, then by Lemma \ref{degree}$(iii)$, either $d_1-1=d_2-1$ or  $d_1-1=d_2-2$ that is, $d_1=d_2$ or  $d_2=d_1+1$  but if $d_2=d_1+1$, then  $n|(d_1+1)^2$, on the other hand  $d_1|n,(d_1+1)|n$  and $gcd(d_1,d_1+1)=1$ which implies that  $d_1(d_1+1)|n$ and these all together  will imply that $d_1|d_1+1$, a contradiction. Hence $d_1=d_2$.
\end{proof}
	An outline of the above proof can be seen in \cite{anderson1999zero} as well. Since $ann(\textbf{x})=ann(\textbf{y})$ if and only if $\textbf{x},\textbf{y}\in \Omega_d$ for some $d\in \Upsilon(n)$, the set $\mathcal{V}=\{\Omega_d: d\in \Upsilon(n)\} $ is precisely the vertex set of $\Gamma_E(\mathbb{Z}_n)$. \\
\begin{lemma}
For the ring $\mathbb{Z}_n$, $\Gamma_E(\mathbb{Z}_n)\cong\Gamma_{Ann}(\mathbb{Z}_n)$.
\end{lemma}
\begin{proof}
	The proof follows from the isomorphism $\psi$ from $ \mathcal{V}$ to $Ann(\mathbb{Z}_n)$ given by $\psi(\Omega_d)=\langle \textbf{d}\rangle$.
\end{proof}

Since all the non-trivial ideals are principal in $\mathbb{Z}_n$ and $\langle \textbf{x} \rangle \langle \textbf{y} \rangle=\langle \textbf{xy} \rangle$ for all $\textbf{x,y}\in \mathbb{Z}_n$. One can attempt to express the structure of the zero-divisor graph of $\mathbb{Z}_n$ as a generalized join of certain elementary simple graphs. The following lemma proves it. 
\begin{lemma}\label{joiniso}
Suppose that $Ann(\mathbb{Z}_n)=\{ \langle {\bf d}_1\rangle,\langle {\bf d}_2\rangle,\cdots,\langle \textbf{d}_{\tau(n)}\rangle\}$, where $d_i$'s are  divisors of $n$. Then $\Gamma(\mathbb{Z}_n )\cong \Gamma_{Ann}(\mathbb{Z}_n)\big[\Lambda_{d_1},\Lambda_{d_2}\cdots,\Lambda_{d_{\tau(n)}}\big]$,\\ where  $\Lambda_{d_i}= \left\{\begin{array}{ll}
	K_{\phi(\frac{n}{d_i})}& ~\textnormal{if} ~n|d_i^2\\
	\overline{K}_{\phi(\frac{n}{d_i})}&~ ~\textnormal{otherwise}.
\end{array} \right.$
\end{lemma}
\begin{proof} 
If the vertex set of  $\Lambda_{d_i}$ is $\{y_i^j: j\in[1,\phi(\frac{n}{d_i})]\}$, then define the vertex set of $\Gamma_{Ann}(\mathbb{Z}_n)\big[\Lambda_{d_1},\cdots ,\Lambda_{d_{\tau(n)}}\big]$ to be  $$U=\bigcup_{i=1}^{\tau(n)} \{(\langle \textbf{d}_i \rangle,y_i^j): y_i^j\in \Lambda_{d_i}, i\in [1,\tau(n)]\}.$$
Now consider a mapping $\psi:V(\Gamma(\mathbb{Z}_n))\rightarrow U$ in such a way that, if $\textbf{x}_j$ is the $j$-th vertex in $ \Omega_{d_i}$ with respect to some order and  if $j \in [1,\phi(\frac{n}{d_i})]$, then define  $\psi(\textbf{x}_j)=(\langle \textbf{d}_i \rangle,y_i^j)$. We claim that 
$\psi$ is a graph isomorphism. Suppose $\textbf{x}$ and $\textbf{y}$ are adjacent in $\Gamma(\mathbb{Z}_n)$.  If $\textbf{x},\textbf{y}\in \Omega_{d_j}$ for some $j$, then $d_i^2|n$ and $\Lambda_{d_i} \cong K_{\phi(\frac{n}{d_i})} $ and hence by definition of generalized join, $\psi(\textbf{x})$ will be adjacent to $\psi(\textbf{y})$. If $\textbf{x}\in \Omega_{d_i}$ and $\textbf{y}\in \Omega_{d_j}$ with $i\ne j$, then $n|d_id_j$ and therefore $\langle \textbf{d}_i \rangle$ and $\langle \textbf{d}_j \rangle$ are adjacent in $\Gamma_{Ann}(\mathbb{Z}_n)$ and again by definition of generalized join of graphs $\psi(\textbf{x})$ is adjacent to $\psi(\textbf{y})$.
Therefore $\psi$ is a graph isomorphism.
\end{proof} 
\begin{lemma}\label{autZn}\cite{anderson1999zero}
 $Aut(\Gamma(\mathbb{Z}_n))=\prod_{d\in \Upsilon(n)} S_{\Omega_d}\cong  \prod_{d\in \Upsilon(n)} S_{\phi(\frac{n}{d})}$.
\end{lemma}

\begin{lemma}\label{orbit}\cite{boutin}
	Let $\mathcal{O}_1,\mathcal{O}_2,\cdots,\mathcal{O}_k$ be the vertex orbits of a graph $\Gamma$ and  $\Lambda_1,\Lambda_2,\cdots, \Lambda_k$ respectively be their  associated	induced subgraphs in $\Gamma$. If  $S_1,S_2, \cdots, S_k$ are the respective determining sets of $\Lambda_1,\Lambda_2,\cdots, \Lambda_k$. Then  $ \bigcup_{i=1}^k S_i$ is a determining set of the  graph $\Gamma$.
\end{lemma}
\begin{theorem}\label{detzn} For the ring $\mathbb{Z}_n$, 
 $Det(\Gamma(\mathbb{Z}_n))=n-\phi(n)-\tau(n)-1$.
\end{theorem}
\begin{proof}
Let $n=4$. Then $Det(\Gamma(\mathbb{Z}_n))=0$, because $\Gamma(\mathbb{Z}_4)=K_1$. In this case,  $n-\phi(n)-\tau(n)-1=4-2-1-1.$ \\
For $n\ne 4$, 
let $\Upsilon(n)=\{d_1,d_2,\cdots,d_{\tau(n)}\}$ and suppose $\Lambda_i$ is the subgraph of $\Gamma(\mathbb{Z}_n)$ induced by  $\Omega_{d_i}$ and let $S_i$ be a  minimum determining set of $\Lambda_i$, that is, $Det(\Lambda_i)=|S_i|$. By Lemma \ref{CIset}, each  $\Lambda_i$ is either a complete graph or an empty graph. In both cases, $S_i$ contains all the elements of $\Omega_{d_i}$, except for exactly one vertex  of it. Hence $|S_i|=|\Omega_{d_i}|-1$.\\
Now, we claim that $S= S_1\cup S_2 \cup \cdots \cup S_{\tau(n)}$ can form a determining set with minimum cardinality for $\Gamma(\mathbb{Z}_n)$. Using Lemma \ref{autZn}, each of $\Omega_{d_i}$ is a vertex orbit of $\Gamma(\mathbb{Z}_n)$. Therefore by Lemma \ref{orbit}, $S$ is a determining set of $\Gamma(\mathbb{Z}_n)$. To prove $S$ is a minimum determining set, let $S'$ be any set of vertices of $\Gamma(\mathbb{Z}_n)$ with $|S'| <|S|$. Then one can find $d_j\in \Upsilon(n)$  and $ \textbf{x},\textbf{v} \in  \Omega_{d_j} $ such that $u,v \notin S'$. Again by Lemma \ref{autZn}, there exists an automorphism of $\Gamma(\mathbb{Z}_n)$, which fixes all the vertices of $S'$ other than $\textbf{x}$ and $\textbf{y}$, but $\textbf{x}$ and $\textbf{y}$ are mapped each other. Thus $S'$ is not a fixing set of $\Gamma(\mathbb{Z}_n)$ and hence by Lemma \ref{fixdet}, $S'$ is not a determining set of $\Gamma(\mathbb{Z}_n)$. It is not hard to verify that $|S|=|\mathcal{Z}(\mathbb{Z}_n)|-|\Upsilon(n)|=n-\phi(n)-\tau(n)-1$.
\end{proof}
\begin{theorem}
	For the ring $\mathbb{Z}_n$, $dim_M(\Gamma(\mathbb{Z}_n))=n-\phi(n)-\tau(n)-1$.
\end{theorem}
\begin{proof} Let $n=p^2$, where $p$ is a prime. Then by Lemma \ref{joiniso}, $\Gamma(\mathbb{Z}_n)$ is a complete graph on $p-1$ vertices. Therefore $dim_M(\Gamma(\mathbb{Z}_n))=p-2$, which is equal to $n-\phi(n)-\tau(n)-1=p^2-(p^2-p)-1-1$.\\
	For any other values of $n$, by Proposition \ref{detmd}, $dim_M(\Gamma(\mathbb{Z}_n))\ge n-\phi(n)-\tau(n)-1$. For the reverse inequality, it is enough to prove that the determining set $S$ defined in the proof of Theorem 1, is a resolving set. Whenever $\textbf{x}$ and $\textbf{y}$ are any two vertices of $\Gamma(\mathbb{Z}_n)$ that are not in $S$, we need to prove $D_S(\textbf{x})\ne D_S(\textbf{y})$. Using definition of $S$, we can find $i\ne j$ such that $\textbf{x} \in \Omega_{d_i}$ and $\textbf{y}\in \Omega_{d_j}$. Then $deg_{\Gamma(\mathbb{Z}_n)}(\textbf{x}) \ne deg_{\Gamma(\mathbb{Z}_n)}(\textbf{y})$. Without loss of generality assume that $deg_{\Gamma(\mathbb{Z}_n)}(\textbf{x}) >deg_{\Gamma(\mathbb{Z}_n)}(\textbf{y})$, then there exist three vertices $\textbf{x},\textbf{y}, \textbf{w}\in \mathcal{Z}(\mathbb{Z}_n)$, $\textbf{x} \ne \textbf{w}$, and $\textbf{y}\ne \textbf{w}$ such that $\textbf{w}\textbf{x}=\textbf{0}$ and $\textbf{w}\textbf{y}\ne \textbf{0}$. Let  $\Omega_d$ be the partite set containing $\textbf{w}$.\\
	If $|\Omega_d| \ne 1$, then there exists $\textbf{w}'\in S\cap \Omega_d$ such that $d(\textbf{x},\textbf{w}')=1$ and $d(\textbf{y},\textbf{w}')\ne 1$  and hence $D_S(\textbf{x})\ne D_S(\textbf{y})$.\\
	Now, if $|\Omega_d|=1$, then by Lemma \ref{Vdsize}(ii), $d=\frac{n}{2}$. Hence $n$ must be even in this case and thus $2$ divides $n$ and $|\Omega_2|>1$. Therefore it is possible to  choose $\textbf{w}'\in S\cap \Omega_2$ such that $d(\textbf{x},\textbf{w}')=2$ and $d(\textbf{y},\textbf{w}')>2$. Therefore $D_S(\textbf{x})\ne D_S(\textbf{y})$ and hence $S$ is a resolving set of $\Gamma(\mathbb{Z}_n)$.
\end{proof}

For example, the zero-divisor graph of $\mathbb{Z}_{315}$ is given below. This graph has  both the determining number and metric dimension  are equal to 160. Note that the grey color connecting $\overline{K}_i$ and  $\overline{K}_j$ represents every vertex of $\overline{K}_i$ is adjacent to every vertex of $\overline{K}_j$.
%-----------------
\begin{figure}[htbp]
  \centering
  \begin{tikzpicture}[scale=0.45]
\draw (0,0)--(0,5);
  \filldraw[line width=4mm, fill = white, draw = gray!30] (0,0) circle (5 cm);
   \draw[line width=4mm,draw = gray!30] (0,0.7)--(0,4.3);
  \draw[line width=4mm,draw = gray!30] (0,-0.7)--(0,-4.3);
  \draw[line width=4mm,draw = gray!30] (3.8,-2)-- (-3.8,-2);
  \draw[line width=4mm,draw = gray!30] (-4.3,-1.4)--(-0.45,4.4) ;
  \draw[line width=4mm,draw = gray!30] (4.3,-1.4)--(0.45,4.4);
  \draw[line width=4mm,draw = gray!30] (3.8,2)-- (-3.8,2);
  \draw[line width=4mm,draw = gray!30]  (-7.5,-4.5)-- (-4.9,-2.5);
  \draw[line width=4mm,draw = gray!30]  (7.5,-4.5)-- (4.9,-2.5);
  \draw[line width=4mm,draw = gray!30]  (0,5.6)--  (0,8.3);
  \filldraw[fill = white] (0,5) circle (0.7 cm);
  \filldraw[fill = white, draw = black] (0,-5) circle (0.7 cm);
  \filldraw[fill = white, draw = black] (4.5,2) circle (0.7 cm);
  \filldraw[fill = white, draw = black] (-4.5,2) circle (0.7 cm);
  \filldraw[fill = white, draw = black] (4.5,-2) circle (0.7 cm);
  \filldraw[fill = white, draw = black] (-4.5,-2) circle (0.7 cm);
    \filldraw[fill = white, draw = black] (-8,-5) circle (0.7 cm);
 \filldraw[fill = white, draw = black] (8,-5) circle (0.7 cm);
  \filldraw[fill = white, draw = black] (0,9) circle (0.7 cm);
  \filldraw[fill = white, draw = black] (0,0) circle (0.7 cm);
   \node at (0,0) {\scriptsize${\overline{K}_{24}}$} ;
  \node at (0,5) {\scriptsize$\overline{K}_{2}$} ;
  \node at  (0,-5){\scriptsize${\overline{K}_{6}}$} ;
  \node at (4.5,2) {\scriptsize${\overline{K}_{12}}$} ;
  \node at (-4.5,2) {\scriptsize${\overline{K}_{8}}$} ;
  \node at (4.5,-2) {\scriptsize${\overline{K}_{4}}$} ;
  \node at  (-8,-5){\scriptsize${\overline{K}_{24}}$} ;
  \node at (8,-5) {\scriptsize${\overline{K}_{36}}$} ;
  \node at (0,9){\scriptsize${\overline{K}_{48}}$} ;
  \node at  (-4.5,-2) {\scriptsize${\overline{K}_{6}}$} ;
  \node at (0,10.5){\small${\Omega_{3}}$} ;
  \node at (0,-6.5){\small${\Omega_{35}}$} ;
  \node at (8,-6.5){\small${\Omega_{5}}$} ;
  \node at (-8,-6.5){\small${\Omega_{7}}$} ;
  \node at (-6,2){\small${\Omega_{21}}$} ;
  \node at (-6,-2){\small${\Omega_{45}}$} ;
  \node at (6,2){\small${\Omega_{15}}$} ;
  \node at (6,-2){\small${\Omega_{63}}$} ;
  \node at (1,1){\small${\Omega_{9}}$} ;
  \node at (1,6){\small${\Omega_{100}}$} ;
\end{tikzpicture}
  \caption{$\Gamma(\mathbb{Z}_{315})$}
  \end{figure}
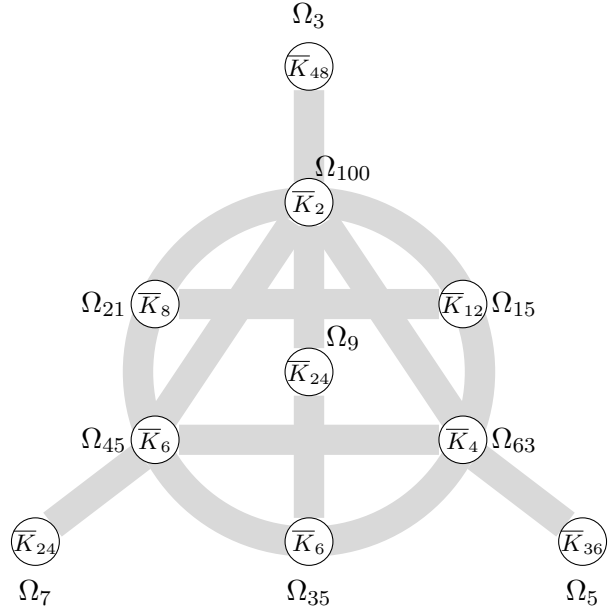
%-----------------

\section{Determining number of zero-divisor graph of  a semisimple ring}
In this section, the determining number and metric dimension of the zero-divisor graph of semisimple  rings are discussed. We use the following notations and results for the same.
\begin{definition}\cite{atiyah}
Let	$R$ be a finite  ring, Then $R$ is said to be semisimple if $R$ is a direct product of a finite number of finite fields.
\end{definition}
Therefore, whenever $R$ is semisimple ring, then $R\cong \prod_{i=1}^k F_i$, where all $F_i$'s are finite field and $k>1$.  Note that if $S$ is a subset of the vertex set of $\Gamma(R)$ of a semisimple ring $R$, then the complement of $S$, denoted by $S^c$, is the set of all vertices in  $\Gamma(R)$, which are not in $S$.
\begin{definition}\label{theta}
\cite{mypaper1} Let $R\cong \prod_{i=1}^k F_i$.
	For any zero-divisor $x=(x_1, \cdots,x_k) $ of $R$, we define a non-empty subset $\Theta_x$ of $[1,k]$ such that $x_i=0$, whenever $i \in \Theta_x$. For any ideal $\mathcal{I}=\prod_{i=1}^kI_i$ of $R$, we define a non-empty subset $\Theta_{\mathcal{I}}$ of  $[1,k]$ such that $I_i= \bf{0}$, whenever $i\in \Theta_{\mathcal{I}}$. Also, we define a subset $\mathcal{I}'= \{x\in \mathcal{I}: \Theta_x=\Theta_{\mathcal{I}}\}$ of $\mathcal{I}$.
\end{definition}

\begin{lemma}\label{genjoinsemi}\cite{mypaper1}	
	Let $R\cong \prod_{i=1}^k F_i$ be a finite semisimple  ring with $Ann(R)=\{\mathcal{I}_1,\mathcal{I}_2,\cdots,\mathcal{I}_m\}$, where $m=2^{k}-2$. For any  $i\in[1,m]$, if $\Lambda_i$ is a totally disconnected graph on $|\mathcal{I}'_i|$ vertices, then  $\Gamma(R)\cong \Gamma_{Ann}(R)\bigl[\Lambda_1,\Lambda_2,\cdots,\Lambda_m \bigr]$.
\end{lemma}
\begin{lemma}\label{autofjoin}
	Let $\Gamma\cong \Lambda[\Lambda_1,\cdots \Lambda_k]$, where $\Lambda$ and $\Lambda_i$'s are any graphs. Then $Aut(\Lambda_i)\times\cdots \times Aut(\Lambda_k) $ can be embedded in $Aut(\Gamma)$.
\end{lemma}
\begin{proof}Suppose that $V(\Lambda)=[1, k]$  and $V(\Gamma)=\{(x_i,y_i^j):y_i^j\in V(\Lambda_i)\}$.
	Consider a mapping $\psi:Aut(\Lambda_1)\times\cdots \times Aut(\Lambda_k) \rightarrow Aut(\Gamma) $ such that $\psi(\sigma_1,\cdots,\sigma_k)(i,y_i^j)=\sigma_i(y_i^j)$. Now if $\psi(\sigma \circ \rho)=\psi(\sigma_1\circ \rho_1,\cdots ,\sigma_k\circ \rho_k)=\mu$, then, $\mu(i,y_i^j)=\sigma_i(\rho_i(y_i^j))$. Thus, $\psi(\sigma \circ \rho )=\psi(\sigma)\circ\psi(\rho)$. Hence $\psi$ is a homomorphism. Also, $ker(\psi)=ker(\sigma_1)\times\cdots \times ker(\sigma_k)=\{e\}$, because all the  automorphisms $\sigma_i$'s are one-one. Hence $\psi$ is an embedding.
\end{proof}

\begin{lemma}\label{lemma1}\cite{mypaper1}
	Let $R\cong \prod_{i=1}^k F_i$ be a semisimple  ring, Then $\Gamma_E(R)\cong \Gamma_{Ann}(R)$.\end{lemma}
 \begin{theorem}\label{detsemi}
	Let $R\cong \prod_{i=1}^{k} F_i$, where $F_i$ is a finite field and $F_i \ncong \mathbb{Z}_2$ for some $i$. Then $Det(\Gamma(R))=|\mathcal{Z}(R)|-2^k+2$.
\end{theorem}
\begin{proof}
From  Lemma \ref{genjoinsemi}, $\Gamma(R)$ is the generalized join of a collection of complete graphs.
Let $S$ be a subset of the vertex set of  $\Gamma(R)$ such that $S$ contains all but one vertex from each of the partite set $\mathcal{I}'$. Since there are $2^k-2$ different partite sets, we have $|S|=|\mathcal{Z}(R)|-2^k+2$. We claim that $S$ is a minimum determining set of $\Gamma(R)$. Suppose $\psi$ is an automorphism of $\Gamma(R)$ that fixes all the vertices of $S$ and let $x \in   S^c$ be such that  $\psi(x) =y$ for some $y\ne x$. Also, assume that $x$ is in $\mathcal{I}'$ of the partition. Now, assume that $|\mathcal{I}'|$=1.  Since there is at least one $j$ such that $F_j\not \cong \mathbb{Z}_2$, there exists  $w\in S$ such that $w$ is adjacent to $x$ but not to $y$. Therefore $\psi(x)$ is not adjacent to $w=\psi(w)$. Hence $\psi $ is not an automorphism of $\Gamma(R)$, a contradiction. Thus $\psi$ fixes all the vertices of $\Gamma(R)$ and $S$ is a determining set of $\Gamma(R)$.\\
Similarly, if $|\mathcal{I}'|>1$, by the procedure that $S$ is defined we must have $y\not \in \mathcal{I}'$. Then  there exists a vertex $w$ either in $S$ or in $\mathcal{I}'$ of the partition with $|\mathcal{I}'|=1$ such that $w$  is adjacent to  $x$, but not to $y$ in $\Gamma(R)$.   
From above, every vertex in $S$ and in  $\mathcal{I}'$ of the partition with $|\mathcal{I}'|=1$  are fixed by $\psi$, therefore $\psi(w)=w$ and hence  $\psi$ is not  an automorphism.
 Considering both the cases together, $\psi(x)=x$ for all $x\in V(\Gamma(R))$. Therefore $S$ is a determining set of $\Gamma(R)$. \\
To prove minimality, suppose  $S_1$ is any subset of the vertex set of $\Gamma(R)$ such that $|S_1|<|\mathcal{Z}(R)|-2^k+2$, then there exist $x,y\in   S_1^c$ such that, $x,y \in \mathcal{J}'$ for some ideal $\mathcal{J}$ of $R$. Then, by Theorem \ref{autofjoin}, there exists an automorphism of $\Gamma(R)$, which fixes all vertices in $S'$; but maps $x$ and $y$ each other. Hence $S_1$ is not a determining set of $\Gamma(R)$. 
\end{proof}
\begin{theorem}
	Let $R\cong \prod_{i=1}^{k} F_i$, where $F_i$ is a finite field and $F_i \ncong \mathbb{Z}_2$ for some $i$. Then $dim_M
	(R)=|\mathcal{Z}(R)|-2^k+2$.
\end{theorem}
	\begin{proof}
Let $S$ be the same as in the proof of Theorem \ref{detsemi}; Since $R$ is not a Boolean ring, $S$ is non-empty. It is enough to prove that, above determining set $S$ is a resolving set as well. Suppose, $x,y$ are distinct vertices from $ S^c$, using the definition of  $S$  we have $\Theta_x\ne \Theta_y$. Then it is possible to find a vertex $w$ in $S$, which is adjacent to one of $\{x,y\}$  but adjacent not to the other one. For, if $\Theta_x\cap  \Theta_y=\{\}$, choose $w$ such that, $\Theta_x^c\subseteq \Theta_w$ and if $\Theta_x\cap  \Theta_y\ne \{\}$, choose $w$ such that, $(\Theta_x\setminus\Theta_y)^c \subseteq \Theta_w$.  Hence, the set $S$  is a resolving set as well, and $dim_M(R)=|\mathcal{Z}(R)|-2^k+2$.
\end{proof}
\begin{theorem}
	Let $R\cong \mathbb{Z}_2^n$ where $n>1$. Then $Det(\Gamma(R))<  \frac{n}{2}+1$.
\end{theorem}
\begin{proof}
First we will consider the cases when $n=2,3$ and $4$ separately.\\
If $n=2$, then $\Gamma(\mathbb{Z}_2^2)\cong K_2$ and hence $Det(\Gamma(\mathbb{Z}_2^2))=1$.\\

	\begin{figure}[!htbp]
	\label{minipage}
\begin{minipage}{6 cm}
	\includegraphics[scale=0.45]{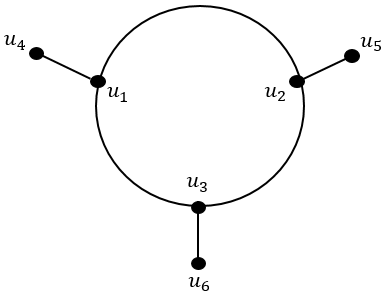}
	
\end{minipage}%
\begin{minipage}{6 cm}
    \includegraphics[scale=0.38]{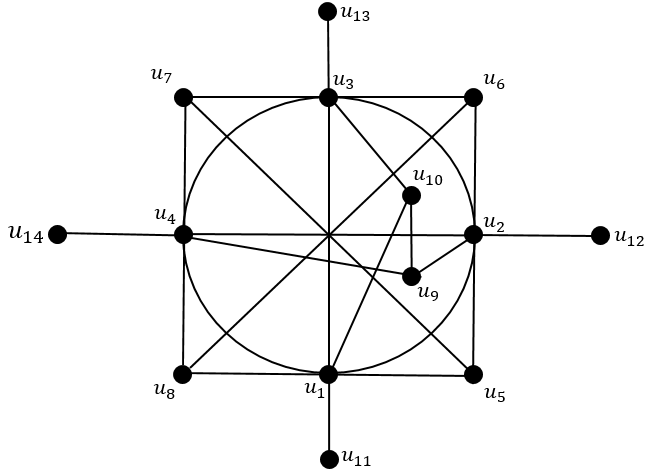}
    
\end{minipage}
	\caption{Zero-divisor graph of $\mathbb{Z}_2^3$ and $\mathbb{Z}_2^4$}
	\end{figure}\vspace{-.2cm}
For $n=3$, we have $\Gamma(\mathbb{Z}_2^3)\cong K_3\odot K_2$, which is given in Figure 2. One can see that   $\{u_1,u_2\}$ is a  minimum determining set of $\Gamma(\mathbb{Z}_2^3)$ and hence the determining number of $\Gamma(\mathbb{Z}_2^4)$ is 2. 
\par For $n=4$, we see that $\{v_6,v_{10}\}$ is a minimum determining set of $\Gamma(\mathbb{Z}_2^4)$ and thus the determining number of $\Gamma(\mathbb{Z}_2^4)$ is $2$.
 \par For $n\ge 5$,  consider the  set $U=\{v\in \mathbb{Z}_2^n : |\Theta_v|=1\}$. Since every vertex of $\Gamma(\mathbb{Z}_2^n)$ is adjacent to at least one of the vertex in $U$,  $U$  itself  is a determining set of $\Gamma(\mathbb{Z}_2^n)$.  In addition, $U$ is precisely the set of the central vertices of $\Gamma(\mathbb{Z}_2^n)$.

Now, for the rest of the proof, we consider two cases: (i) for even $n=2t, t\geq 3$ and (ii) for odd $n=2s+1, s\geq 2$.\\ 
Consider the portion of the graph $\Gamma(\mathbb{Z}_2^{2t})$ illustrated in Figure 3. Let $U=\{v_1,v_2,\cdots, v_{2t}\}$ be the set of central vertices of $\Gamma(\mathbb{Z}_2^{2t})$. Then, choose the vertices $u_1,u_2,\cdots u_{t}$ such that for each $i\in[1,t]$, the vertex $u_i$ is adjacent to exactly three of the central vertices $v_{2i-1},v_{2i}$ and $v_{2i+1}$ as given in Figure $3$. Note that it is considered $v_{2k+1} $ is the same as that of $v_1$. Let $S=\{u_1,u_2,\cdots u_t\}$ be set of the vertices of $\Gamma(\mathbb{Z}_2^{2k} )$, as shown in Figure 3. Suppose an automorphism of $\Gamma(\mathbb{Z}_2^{2k})$ fixes all of the vertices of $S$, then it forces all the central vertices to be fixed. Hence $S$ constitutes a determining set of $\Gamma(\mathbb{Z}_2^{2k})$. 
	\begin{figure}[!htbp]
	\label{minipage}
\begin{minipage}{6 cm}
\begin{tikzpicture}[scale=0.43]
	\draw (-7,-7) rectangle ++(14,14);
	\tkzLabelPoint[below,xshift=10mm,yshift=10](135:11){$\Gamma(\mathbb{Z}_2^{2k})$}
	\filldraw[fill=lightgray,draw=none](0,0) circle (4);
%	\draw[dashed] (120:2.5)--(120:-2.5);
%	\draw[dashed] (30:2.5)--(30:-2.5);
	\draw [dotted,domain=-70:-20](0,0) plot ({4*cos(\x)}, {4*sin(\x)});
	\draw [dotted,domain=-150:-200](0,0) plot ({4*cos(\x)}, {4*sin(\x)});
	\draw [dotted,domain=-20:-70] plot ({6*cos(\x)}, {6*sin(\x)});
	\draw [dotted,domain=-150:-200] plot ({6*cos(\x)}, {6*sin(\x)});
	\draw [black,domain=0:135] plot ({4*cos(\x)}, {4*sin(\x)});
	\draw [black,domain=225:270] plot ({4*cos(\x)}, {4*sin(\x)});
%	\draw [black,domain=225:275] plot ({4*cos(\x)}, {4*sin(\x)});
%	\draw [dashed,domain=130:330] plot ({6*cos(\x)}, {6*sin(\x)});
	\fill [black] (0,0) +(90:4) circle (4pt);
	\fill [black] (0,0) +(67.5:4) circle (4pt);
	\fill [black] (0,0) +(112.5:4) circle (4pt);
	\fill [black] (0,0) +(135:4) circle (4pt);
	\fill [black] (0,0) +(0:4) circle (4pt);
	\fill [black] (0,0) +(22.5:4) circle (4pt);
	\fill [black] (0,0) +(45:4) circle (4pt);
	\fill [black] (0,0) +(67.5:6) circle (4pt);
	\fill [black] (0,0) +(22.5:6) circle (4pt);
	\fill [black] (0,0) +(112.5:6) circle (4pt);
	\fill [black] (0,0) +(247.5:6) circle (4pt);
	\fill [black] (0,0) +(247.5:4) circle (4pt);
	\fill [black] (0,0) +(270:4) circle (4pt);
	\fill [black] (0,0) +(225:4) circle (4pt);
 	\draw (135:4cm)--(112.5:6cm)--(90:4cm)--(67.5:6cm)--(45:4cm)--(22.5:6cm)--(0:4cm);
 	\draw (270:4cm)--(247.5:6)--(225:4);
 	\draw (247.5:4cm)--(247.5:6cm);
 	\draw (67.5:4)--(67.5:6);
 	\draw (22.5:4)--(22.5:6);
 	\draw (112.5:4)--(112.5:6);
 		\tkzLabelPoint[below](90:4){$v_1$};
 		\tkzLabelPoint[below,xshift=-1mm](67.5:4){$v_2$}
 		\tkzLabelPoint[left,xshift=0mm,yshift=-2mm](45:4){$v_3$}
 		\tkzLabelPoint[left,xshift=0mm,yshift=0mm](22.5:4){$v_4$}
 		\tkzLabelPoint[left,xshift=0mm,yshift=0mm](0:4){$v_5$}
 		\tkzLabelPoint[above,xshift=.5mm,yshift=0mm](270:4){$v_{2i-1}$}
 		\tkzLabelPoint[above,xshift=0mm,yshift=0mm](247.5:4){$v_{2i}$}
 		\tkzLabelPoint[above,xshift=1.6mm,yshift=0.15mm](225:4){$v_{2i+1}$}
 		\tkzLabelPoint[below,xshift=.05mm,yshift=0mm](112.5:4){$v_{2k}$}
\tkzLabelPoint[below,xshift=4mm,yshift=0mm](135:4){$v_{2k-1}$}
 		%\tkzLabelPoint[below,xshift=4mm,yshift=0mm](-45:6.4){$S$}
 		\tkzLabelPoint[above,xshift=1mm,yshift=1mm](67.5:6){$u_1$}
 		\tkzLabelPoint[above,xshift=4mm,yshift=-.1mm](22.5:6){$u_2$} 		\tkzLabelPoint[below,xshift=-.3mm,yshift=-2mm](245.5:6){$u_i$}
 		\tkzLabelPoint[above,xshift=-2.5mm,yshift=1.2mm](112.5:6){$u_k$}
		\end{tikzpicture}
\end{minipage}%
\begin{minipage}{6 cm}
   \begin{tikzpicture}[scale=0.43]
		\draw [draw=black] (-7,-7) rectangle ++(14,14);
	\draw [dotted,domain=-20:-70] plot ({6*cos(\x)}, {6*sin(\x)});
		\tkzLabelPoint[below,xshift=10mm,yshift=10](135:11){$\Gamma(\mathbb{Z}_2^{2k+1})$}
	\draw [dotted,domain=-150:-200] plot ({6*cos(\x)}, {6*sin(\x)});
		\draw [dotted,domain=-150:-190] plot ({4*cos(\x)}, {4*sin(\x)});
			\draw [dotted,domain=-70:-20] plot ({4*cos(\x)}, {4*sin(\x)});
		\fill[fill=lightgray,draw=none](0,0) circle (4);
		\draw [black,domain=0:157.5] plot ({4*cos(\x)}, {4*sin(\x)});
		\draw [black,domain=225:270] plot ({4*cos(\x)}, {4*sin(\x)});
	%	\draw [dashed,domain=150:330] plot ({6*cos(\x)}, {6*sin(\x)});
		\fill [black] (0,0) +(90:4) circle (4pt);
		\fill [black] (0,0) +(67.5:4) circle (4pt);
		\fill [black] (0,0) +(112.5:4) circle (4pt);
		\fill [black] (0,0) +(135:4) circle (4pt);
		\fill [black] (0,0) +(157.5:4) circle (4pt);
%		\fill [black] (0,0) +(180:4) circle (4pt);
		\fill [black] (0,0) +(0:4) circle (4pt);
		\fill [black] (0,0) +(22.5:4) circle (4pt);
		\fill [black] (0,0) +(45:4) circle (4pt);
		\fill [black] (0,0) +(67.5:6) circle (4pt);
		\fill [black] (0,0) +(22.5:6) circle (4pt);
		\fill [black] (0,0) +(135:6) circle (4pt);
				\fill [black] (0,0) +(270:4) circle (4pt);
				\fill [black] (0,0) +(245.5:4) circle (4pt);
				\fill [black] (0,0) +(225:4) circle (4pt);
				\fill [black] (0,0) +(245.5:6) circle (4pt);
				
				%\tkzLabelPoint[below,xshift=4mm,yshift=0mm](-45:6.4){$S$}
		\draw (90:4cm)--(67.5:6cm)--(45:4cm)--(22.5:6cm)--(0:4cm);
		\draw (270:4cm)--(245.5:6cm)--(225:4cm);
		\draw(245.5:4cm)--(245.5:6cm);
		\draw(67.5:4)--(67.5:6);
		\draw(22.5:4)--(22.5:6);
		\draw(135:4)--(135:6);
		\draw (157.5:4)--(135:6)--(112.5:4);
		\tkzLabelPoint[below](90:4){$v_1$};
		\tkzLabelPoint[below,xshift=0mm](67.5:4){$v_2$}
		\tkzLabelPoint[below,xshift=-2mm](45:4){$v_3$}
		\tkzLabelPoint[below,xshift=-2mm,yshift=0mm](22.5:4){$v_4$}
		\tkzLabelPoint[left,xshift=0mm,yshift=-.5mm](0:4){$v_5$}
		\tkzLabelPoint[below,xshift=1mm,yshift=0mm](112.5:4){$v_{2k+1}$}
		\tkzLabelPoint[below,xshift=2mm,yshift=0mm](135:4){$v_{2k}$}
		\tkzLabelPoint[right,xshift=0mm,yshift=-.5mm](157.5:4){$v_{2k-1}$}
		\tkzLabelPoint[above,xshift=1.2mm,yshift=.5mm](67.5:6){$u_1$}
		\tkzLabelPoint[above,xshift=3.5mm,yshift=-.8mm](22.5:6){$u_2$}
		\tkzLabelPoint[above,xshift=-3mm,yshift=.5mm](135:6){$u_{k}$}
		\tkzLabelPoint[below,xshift=-1mm,yshift=-2.5mm](245.5:6){$u_{i}$}
		\tkzLabelPoint[above,xshift=0mm,yshift=0mm](245.5:4){$v_{2i}$}
		\tkzLabelPoint[above,xshift=1mm,yshift=0mm](270:4){$v_{2i-1}$}
		\tkzLabelPoint[above,xshift=-1mm,yshift=0mm](225:4){$v_{2i+1}$}
		\end{tikzpicture}
    
\end{minipage}
	\caption{ $\Gamma(\mathbb{Z}_2^{n})$  and determining set $\{u_1,u_2,\cdots,u_k\}$}
	\end{figure}
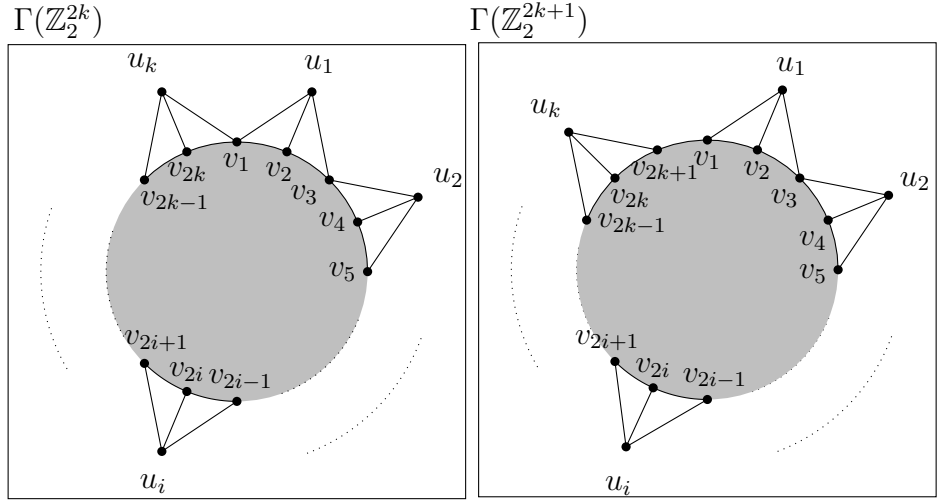

Similarly, for the case $n=2k+1$, Figure 3 illustrates a portion of   $\Gamma(\mathbb{Z}_2^{2k+1})$. Choose a set of vertices $S=\{u_1,u_2,\cdots ,u_k\}$ of $\Gamma(R)$ as shown in Figure 4, with the property that $u_i$ is adjacent to exactly three of the central vertices, namely $v_{2i-1},v_{2i}$ and $v_{2i+1}$ for $i=1,2,\cdots,k$. Then fixing all the vertices of $S$ forces to fix all the central vertices and hence fixes all the vertices of $\Gamma(\mathbb{Z}_2^{2k+1})$.  

Therefore, $Det(\Gamma(\mathbb{Z}_2^n))=\lfloor \frac{n}{2} \rfloor< \frac{n}{2}+1$,  $n\ge 2$.
\end{proof} 

Therefore, $Det(\Gamma(R))=dim_M(\Gamma(R))$ if $R$ is either $\mathbb{Z}_n$ or any non-Boolean reduced ring. A natural question arises is: does there exist a  ring $R$ such that  $Det(\Gamma(R))\ne dim_M(\Gamma(R))$? We have an affirmative answer for this, for the ring $\mathbb{Z}_2^5$,  $Det(\Gamma(\mathbb{Z}_2^5))\le 2$; but, from \cite{rajapirzada}, $dim_M(\Gamma(\mathbb{Z}_2^5))=5$. Therefore, $Det(\Gamma(\mathbb{Z}_2^5))\ne dim_M(\Gamma(\mathbb{Z}_2^5))$.
The following problem is still interesting to study further.
\begin{problem}
Classify the structural properties of  a rings $R$ for which $Det(\Gamma(R))\ne dim_M(\Gamma(R))$.
\end{problem}

\section{Size of the determining number }
Another interesting discussion in this topic is the comparison of the size, that is, the finiteness of a determining number $Det(\Gamma(R))$, with respect to the size of the given ring $R$.
In \cite{pirzu}, Pirzada et al. establish an elegant result that the metric dimension of the zero-divisor graph  of  an arbitrary commutative ring is finite if and only if the ring itself is finite. In the case of determining number, a detailed investigation is required. 
\begin{problem}  
Does there exist an infinite  ring $R$, which is not an integral domain, such that $Det(\Gamma(R)) < \infty$?
\end{problem}
In 2011, Spiroff and Wickam\cite{spiroff2011zero}  associated the vertices of a compressed zero-divisor graph and  associated prime ideals of  a  Noetherian ring. In  \cite{spiroff2011zero}, there  is a natural injective
map from $Ass(R)$ to the vertex set of $\Gamma_E(R)$ given by $\psi(\mathfrak{P})=[y]$, where $\mathfrak{P} = ann(y)$. The following theorem guarantees the size of the determining number of certain zero-divisor graphs, is not finite. 
\begin{theorem}
Let $R$ be any infinite  ring, which is not an integral domain. If $\Gamma_E(R)$ is a finite, then $Det(\Gamma(R))$ is not finite.
\end{theorem}
\begin{proof}
	If  $R$ is an infinite  ring and not an integral domain, then by Lemma \ref{Ganes}, $\mathcal{Z}(R)$ is also an infinite set. If $V(\Gamma_E(R))$ is a finite graph, then at least one of the equivalence class $[x]$ contains an infinite number of zero-divisors. Note that for any two distinct vertices $u$ and $v$ in $[x]$,  $N(u)\setminus\{v\}=N(v)\setminus\{u\}$ in  $\Gamma(R)$. Therefore, any determining set of $\Gamma_E(R)$  includes all but one vertex from $[x]$. Thus, all but a finite number of vertices of $V(\Gamma(R))$  are included in any of the determining sets of $\Gamma(R)$. Hence, $Det(\Gamma(R))$ is not  a  finite number.
\end{proof}
%\begin{theorem}
	%Let $R$ be an infinite Noetherian  ring and there is no vertex $[x]$ in $\Gamma_E(R)$ such that $deg_{\Gamma(\mathbb{Z}_n}[x]=\infty$ and $\psi^{-1}([x])$ is a associated prime maximal in $Ann(R)$. Then $Det(\Gamma(R))=\infty$.
%\end{theorem} 
%\begin{proof}
%Using the Proposition 2.2 of \cite{spiroff2011zero}, $\Gamma_E(R)$ is finite.
%\end{proof}
The converse of the above theorem need not be true.
The determining number of the zero-divisor graph of an infinite ring $R$ need not be finite, even if $\Gamma_E(R)$ is not a finite graph. For example, consider $R=\mathbb{Z}_2[x,y,z]/(x^2,y^2)$  is an infinite ring and $\Gamma_E(R)$ is  infinite, but $Det(\Gamma(R))$ is not finite.

\begin{theorem}\cite{mulay}\label{mulayclass}
Let $R$ be a  Noetherian ring and suppose all the equivalence classes $[x]$ with respect to the relation $\sim$  has finite cardinality. Then $R$ is a finite ring.
\end{theorem}
The following theorem will settle the above problem when the ring is Noetherian and the problem is still open for the non-Notherian rings.
\begin{theorem}
Let $R$ be Noetherian. Then $Det(\Gamma(R))$ is finite if and only if $R$ is finite.
\end{theorem}
\begin{proof}
By Theorem \ref{mulayclass}, if $R$ is   infinite, then at least one of the  equivalence class $[x]$ has infinite elements. Therefore, any determining set of $\Gamma(R)$ contains all but  one vertex of that class $[x]$. Hence $Det(\Gamma(R))$ is not a finite number.
\end{proof}
\section{Miscellaneous results}
In this section, we discuss certain interesting results in connection  with determining number and metric dimension of general graphs.
When we write $\Gamma$ is isomorphic to  $\Lambda \bigl[\Lambda_{1},\Lambda_{2},\cdots ,\Lambda_{k} \bigr]$, we mean $\Lambda$ is a finite graph with the vertices $1,2,\cdots, k$ and $\Lambda_{1},\Lambda_{2},\cdots, \Lambda_{k}$ are  finite graphs. Also, if $S$ is any subset of $V(\Lambda_i)$, then we use $S^{\dagger}$ to denote the corresponding subset $\{(i,v):v\in S\}$ of $V(\Gamma)$. We also note that the number of elements in $S$ and $S^{\dagger}$ are the same.  
\begin{theorem}
Let  $\Gamma \cong \Lambda \bigl[\Lambda_{1},\Lambda_{2},\cdots , \Lambda_{k} \bigr]$, where $\Lambda_i$ is either a complete graph or an empty graph and  $deg_{\Gamma}((i,u))\ne deg_{\Gamma}((j,v))$  whenever $i\ne j$.  Then $Det(\Gamma)=(\sum_{i=1}^{k} |V(\Lambda_{i})^\dagger|)-k$.
\end{theorem}
\begin{proof}
 As per the notation defined above, we have  $V(\Gamma)=\cup_{i=1}^kV(\Lambda_i)^\dagger$. Clearly by the hypothesis, $ \{V(\Lambda_i)^\dagger\}_{i=1}^k$        are the vertex orbits of $\Gamma$.  Let $S$ be a subset of $V(\Gamma)$ such that for each $i$, the  cardinality of the set $S\cap V(\Lambda_i)^{\dagger}$ is equal to one less than the cardinality of the set $V(\Lambda_i)$ that is, $|S\cap V(\Lambda_i)^{\dagger}|=|V(\Lambda_i)|-1$. Since $|S\cap V(\Lambda_i)^{\dagger}|=|V(\Lambda_i)|-1$, we have $S\cap V(\Lambda_i)$ is a determining set of the subgraph of $\Gamma$ induced by $V(\Lambda_i)^{\dagger}$. Therefore by Lemma \ref{orbit}, $S=\cup_{i=1}^k(S\cap V(\Lambda_i)^{\dagger})$ is a determining set of $\Gamma$. Now, we claim that $S$ is a minimum determining set of $\Gamma$. Let  $S'$ be any determining set of $\Gamma$ with  smaller cardinality than $S$. Then there exists $j\in [1,k]$ such that  $|V(\Lambda_j)^{\dagger}\setminus S'|\ge2$. Let $(j,u)$ and $(j,v)$ be two vertices in $V(\Lambda_j)^{\dagger}\setminus S'$. Then by Theorem \ref{autofjoin}, there exists an automorphism, which fixes all the vertices except $(j,u)$ and $(j,v)$. Thus $S'$ is not a determining set of $\Gamma$ and hence $S$ is a minimum determining set of $\Gamma$. Therefore $Det(\Gamma)=(\sum_{i=1}^{k} |V(\Lambda_{i})|)-k$
\end{proof}
	\begin{theorem}
		Suppose that $\Gamma \cong \Lambda \bigl[\Lambda_{1},\Lambda_{2},\cdots , \Lambda_{k} \bigr]$, where $\Lambda_i$ is a vertex-transitive  graph.  Then $Det(\Gamma)=\sum_{i=1}^{k}Det(\Lambda_{i})$.
	\end{theorem}
	\begin{proof}
	 Since  $\Lambda_{i}$ is a  vertex-transitive graph, the set, $V(\Lambda_{i})^{\dagger}=\{(i,v):v\in V(\Lambda_{i})\} $ will be a vertex orbit of the graph $\Gamma$. Suppose $S_i$ is a minimum determining set of $\Lambda_{i}$, that is, $|S_i|=Det(\Lambda_{i})$. Then the set $S_i^{\dagger}$ is a minimum determining set of the subgraph of $\Gamma$ induced by $V(\Lambda_i)^{\dagger}$. Now using Lemma \ref{orbit}, $S=\bigcup_{i=1}^kS_i^{\dagger}$ is a determining set of $\Gamma$.\\
	 To prove that $S$ is a minimum determining set, assume that $S'$ is a determining set of $\Gamma$ such that $|S'|<|S|$. Then there exists a vertex  $j$ of $\Lambda$ such that $|V(\Lambda_{j})^{\dagger}\cap S'|<|S_i^{\dagger}| $. Then the corresponding set $\{v\in V(\Lambda_j): (j,v)\in V(\Lambda_{j})^{\dagger}\cap S'\}$ forms a determining set of $\Lambda_i$. This contradicts the fact $Det(\Lambda_{j})=|S_j|$. Hence, $Det(\Gamma)=|S|=\sum_{i=1}^{k} Det(\Lambda_{i})$.
	\end{proof}
 We already found that for any graph $\Gamma$, $Det(\Gamma)\le dim_M(\Gamma)$. It is quite fair to ask, how large can  this difference be. In \cite{boutin} Boutin posed   the problem : 
 \begin{problem}\cite{boutin}
 Can the difference between the determining number of a graph and the size of a smallest  resolving set of a graph be arbitrarily large?
 \end{problem}
 We give an affirmative answer to the above problem.
 \begin{center}
 	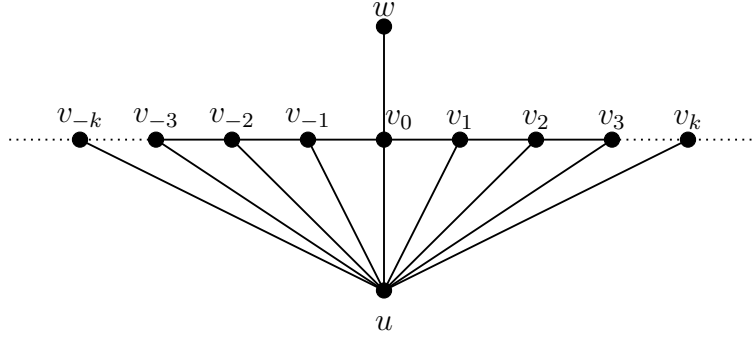
\begin{figure}[!htbp]
 		\begin{tikzpicture}
 		\draw[line width=0.25mm] (-3,0)-- (3,0);
 		\draw[line width=0.25mm,dotted] (3,0)-- (5,0);
 		\draw[line width=0.25mm,dotted] (-3,0)-- (-5,0);
 	\draw[line width=0.25mm] (0,-2)-- (0,1.5); 
 	\draw[line width=0.25mm] (0,-2)-- (1,0); 
 		\draw[line width=0.25mm] (0,-2)-- (2,0); 
 			\draw[line width=0.25mm] (0,-2)-- (3,0); 
 			\draw[line width=0.25mm] (0,-2)-- (4,0);
 			\draw[line width=0.25mm] (0,-2)-- (-1,0);
 			\draw[line width=0.25mm] (0,-2)-- (-2,0);
 			\draw[line width=0.25mm] (0,-2)-- (-3,0);
 			\draw[line width=0.25mm] (0,-2)-- (-4,0);
 			\filldraw[fill = black, draw = black] (-4,0) circle (0.1 cm);
 			\filldraw[fill = black, draw = black] (-3,0) circle (0.1 cm);
 			\filldraw[fill = black, draw = black] (-2,0) circle (0.1 cm);
 			\filldraw[fill = black, draw = black] (-1,0) circle (0.1 cm);
 			\filldraw[fill = black, draw = black] (0,0) circle (0.1 cm);
 			\filldraw[fill = black, draw = black] (1,0) circle (0.1 cm);
 			\filldraw[fill = black, draw = black] (2,0) circle (0.1 cm);
 			\filldraw[fill = black, draw = black] (3,0) circle (0.1 cm);
 			\filldraw[fill = black, draw = black] (4,0) circle (0.1 cm);
 			\filldraw[fill = black, draw = black] (0,1.5) circle (0.1 cm);
 			\filldraw[fill = black, draw = black] (0,-2) circle (0.1 cm);

 			\tkzDefPoint(0,0){0}
 			\tkzDefPoint(1,0){1}
 			\tkzDefPoint(2,0){2}
 			\tkzDefPoint(3,0){3}
 			\tkzDefPoint(4,0){4}
 			\tkzDefPoint(-1,0){-1}
 			\tkzDefPoint(-2,0){-2}
 			\tkzDefPoint(-3,0){-3}
 			\tkzDefPoint(-4,0){-4}
 			\tkzDefPoint(0,1.5){w}
 			\tkzDefPoint(0,-1.5){u}
 			\tkzLabelPoint[above](1){$v_1$}
 			\tkzLabelPoint[above](2){$v_2$}
 			\tkzLabelPoint[above](3){$v_3$} 
 			\tkzLabelPoint[above](4){$v_k$} 
 			\tkzLabelPoint[above,xshift=2mm](0){$v_0$} 
 			\tkzLabelPoint[above](-1){$v_{-1}$} 
 			\tkzLabelPoint[above](-2){$v_{-2}$} 
 			\tkzLabelPoint[above](-3){$v_{-3}$} 
 			\tkzLabelPoint[above](-4){$v_{-k}$} 
 			\tkzLabelPoint[above](w){$w$} 
 			\tkzLabelPoint[below,yshift=-7mm](u){$u$} 
 				
 		\end{tikzpicture}\caption{$\Gamma$ with $Det(\Gamma)=1$ and $dim_M(\Gamma)$ is not finite}
 	\end{figure}\vspace{-.2cm}
 \end{center}

The above graph $\Gamma$ has the automorphism group isomorphic to $C_2$, the cyclic group with two elements. Therefore, $Det(\Gamma)=1$. Now, if $S$ is any resolving set of $\Gamma$ with a finite number of vertices, then one can find a positive integer $n$ such that $v_i\notin S$ whenever $|i|\ge n$. Therefore, $D_S(v_{-k})=D_S(v_{k})$ with $k=n+1$ and hence, $S$ not  a resolving set. Thus, $dim_M(\Gamma)$ can not  be finite. \\

The above example guarantees that the determining number of an infinite graph can be finite. One may think that for an infinite graph $\Gamma$, the finiteness of $Aut(\Gamma)$ enforces the determining number of $\Gamma$ to be finite. The following example illustrates that there may exist an infinite graph with an infinite automorphism group having finite determining number.\\

Consider a graph $\Lambda$ with $V(\Lambda)=\mathbb{N}$ and $ E(\Lambda)=\{(n,n+1):n\in \mathbb{N}\}$.
Note that $\Lambda$ does not have any non-trivial automorphism.
Now we can take infinitely many copies of this graph, namely the graph $\Lambda_1=(V_1,E_1)$ with $V_1=\mathbb{N}\times\mathbb{N},E_1=\{((i,j),(i,j+1)):i,j\in \mathbb{N}\}$.
As we do not have any non-trivial automorphisms inside the connected components, we have that the automorphisms are exactly the permutations that permute the connected components  while keeping each of these components unaltered, and thus $Aut(\Lambda_1)\cong S_{\mathbb{N}}$.\\

Now, add the vertices $V_2=\{v_{\sigma}:\sigma\in S_{\mathbb{N}}\}$ and $E_2=\{(v_{\sigma},(i,j)):v_{\sigma}\in V_2,j\le \sigma^{-1}(i)\}$.
Basically, what we do here is, add a new vertex  corresponding to each permutation on $\mathbb{N}$, and for each such vertex $v_{\sigma}$, we order the connected components by $\sigma(1),\sigma(2),...$ and this vertex $v_{\sigma}$ is made adjacent to the $i-$th components (after ordering with respect to $\sigma$) with $i$ new edges, where the edges go to the first $i$ vertices in the component.\\

Again define $V_3$ and $E_3$ such that $V_3=\{u,u^{'}\}, E_3=\{(u,u^{'})\}\cup\{(u,v):v\in V_2\}$. 
Thus, we define $\Gamma$ with the vertex set  $V= \cup_{i=1}^3 V_i$ and the edge set $E=\cup_{i=1}^3 E_i$. 
 The step-by-step definition of the graph $\Gamma$  is summarized in the following table.\\

\begin{tabular}{c} 
\centering
 \begin{tabularx}{12cm}{X}\hline
Construction of $\Gamma$ with the set of vertices and the set of edges \\  \hline
stage 1. $V_1=\mathbb{N}\times \mathbb{N}$, $E_1=\{((i,j),(i,j+1)):i,j\in \mathbb{N}\}$ \\ \hline
stage 2.  $\cup_{i=1}^2 V_i$, $\cup_{i=1}^2 E_i$, where $V_2=\{v_{\sigma}:\sigma\in S_{\mathbb{N}}\}$ and $E_2=\{(v_{\sigma},(i,j)):v_{\sigma}\in V_2,j\le \sigma^{-1}(i)\}$ \\ \hline
stage 3. $\cup_{i=1}^3 V_i$, $\cup_{i=1}^3 E_i$,  where $V_3=\{u,u'\}$, $E_3=\{(u,u^{'})\}\cup\{(u,v):v\in V_2\}$  \\ \hline
 \end{tabularx}
\end{tabular}
\vskip .2cm \noindent
Now, let $\psi$ be an automorphism of $\Gamma$. Since $u'$ is the only vertex of degree 1, $\psi(u^{'})=u^{'}$. Also, we must have $\psi(u)=u$.\\
As $N(u)=N(\psi(u))=V_2$, we must have  $\psi(V_2)=V_2$. Consequently,
$w\in V_2$ implies $\psi(w)\in V_2$ and therefore, $w\in V_1$ if and only if $\psi(w)\in V_1.$

It is to be noted that $\psi$ again permutes connected components in $V_1$. On contrary, if $\psi(i,1)=(i',j)$ with $j>1$, then $\psi(i,1)$ has two neighbors in $V_1$. Therefor,  $(i,1)$ must have  2 neighbors in $V_1$, which is not true, so $\psi(i,1)=\psi(i',1)$ for some $i'\in \mathbb{N}$ and it easily follows by induction that $\psi(i,j)=\psi(i',j')\implies j=j'$. Therefore every automorphism of $\Gamma$ permutes the  connected components of $V_1$, and hence the group of automorphisms of $\Gamma$ restricted to the domain $V_1$  is isomorphic to $S_{\mathbb{N}}$.\\

Let $\sigma \in S_\mathbb{N}$, which maps $\Gamma_i$ to $\Gamma_{\sigma(i)}$, where $\Gamma_i$ is a component  of $\Lambda_1$, which is a subgraph of $\Gamma$ by taking only the vertices $\{(i,j):j\in \mathbb{N}\}$, and $\psi$ be an automorphism of $\Gamma$ whose restriction to $V_1$ is $\sigma$, we claim that $\Gamma$ is unique.\\
As the number of edges from $v_{\rho}$ of $ V_2$ to $\Gamma_i$  must be the same as the number of edges between $\psi(v_{\rho})$ and $\Gamma_{\sigma(i)}$.  Since there are $\rho^{-1}(i)$ number of edges from $v_{tau}$ to $\Gamma_i$, we have no choice, except to set $\psi(v_\rho)=v_{\sigma \circ \rho}$.\\
 Therefore, $Aut(\Gamma)\cong S_{\mathbb{N}}$, and notice that the automorphism, $\psi$ corresponding to $\sigma\in S_{\mathbb{N}}$ maps $v_{id}$ to $v_{\sigma}$. That is, every automorphism of $\Gamma$, can be represented uniquely by a $\sigma\in S_{\mathbb{N}}$, and the vertex $v_{id}\in V_2\subset V$ has the unique image $v_{\sigma}$.
Thus, by the definition of determining set, $S=\{v_{id}\}$ is a determining set of $\Gamma$.\\ Therefore, a graph with an automorphism group of infinite order can have a finite determining number.
The following problem will be an interesting one to work with.
\begin{problem}
Determine necessary conditions for an infinite graph containing a  determining set with a finite number of vertices?
\end{problem}
\section{Conclusion and Scope}
In this article, we adopt a different approach to understanding the structure of zero-divisor graphs. The larger  $\Gamma(\mathbb{Z}_n)$  is expressed as a generalized join of smaller $\Gamma_{Ann}(\mathbb{Z}_n)$ with a collection of complete and empty graphs. By using this structural identification, the determining number and metric dimension of zero-divisor graph of $\mathbb{Z}_n$, is determined and found both to be equal. 

Similarly, the determining number and metric dimension of zero-divisor of finite non-Boolean semisimple rings, and Boolean rings are predicted. 

Also, it is recognized that there are rings such that determining number and metric dimension of corresponding zero-divisor graphs are different numbers. The question: does there exist an infinite zero-divisor graph with a finite determining number? is taken up and is partially settled for Noetherian rings.  Finally, we emphasize on problems related to the determining number of general graphs aside from zero-divisor graphs. The open problem by 
Boutin\cite{boutin} is settled. A  graph having an infinite automorphism group with a finite determining number is also found. There is a scope for further investigation such as deciding if any  infinite zero-divisor graph has a finite determining number or identifying the different graph structures, which govern the size of determining number.

\end{document}